\documentclass[a4paper,11pt]{amsart}

\usepackage[square,sort,comma,numbers]{natbib}

\usepackage{mathptmx}
\usepackage{hyphenat}

\usepackage{rotating}
\usepackage{floatpag}
\rotfloatpagestyle{empty}

\usepackage{amssymb,amsfonts,amsmath,amsthm,amscd,mathrsfs}
\usepackage{psfrag}

\usepackage{hyperref}

\theoremstyle{plain}

\newtheorem*{problem*}{Problem}

\newtheorem{prop}{Proposition}[section]
\newtheorem{PROB}[prop]{Problem}   
\newtheorem*{theorem*}{Theorem}

\theoremstyle{definition}

\usepackage{latexsym}
\usepackage{xspace}
\usepackage{enumerate, paralist}

\DeclareMathOperator{\Aut}{\mathrm{Aut}}
\DeclareMathOperator{\Inn}{\mathrm{Inn}}

\DeclareMathOperator{\QQ}{\mathbf{Q}}

\newcommand{\tdlc}{tdlc\@\xspace}

\numberwithin{equation}{section}

\begin{document}
\title[Future directions in locally compact groups]{Future directions in locally compact groups:\\ a tentative problem list}

\author{{P.-E. Caprace and N. Monod}}
\maketitle



\begin{flushright}
\begin{minipage}[t]{0.7\linewidth}\itshape\small
She went on saying to herself, in a dreamy sort of way,``Do cats eat bats? Do cats eat bats?'' and sometimes, ``Do bats eat cats?'' for,  you see, as she couldn't answer either question, it didn't much matter which way she put it. 

\hfill\upshape (Lewis Carroll, \emph{Alice's Adventures in Wonderland}, 1865.)
\end{minipage}
\end{flushright}

\medskip

This text is the preprint version of the concluding chapter for the book \emph{New Directions in Locally Compact Groups}~\cite{book} published by Cambridge University Press in the series \emph{Lecture Notes of the LMS}. The recent progress on locally compact groups surveyed in that volume also reveals the considerable extent of the unexplored territories. Therefore, we wish to conclude it by mentioning a few open problems related to the material covered in the book and that we consider important at the time of this writing.

We shall group problems along the themes indicated in the following table of contents; each problem is briefly discussed and accompanied by a list of relevant references for further reading.

\tableofcontents

\section{Chabauty limits}
Recall that the collection $\mathbf{Sub}(G)$ of all closed subgroups of a locally compact group $G$ carries a natural topology, the  \textbf{Chabauty topology}\index{Chabauty space}, for which it is a compact space.

\begin{PROB}\label{prob:LimitAmenable}
Let $G$ be a locally compact group. Is the collection of closed amenable subgroups of $G$ a closed subset of the Chabauty space $\mathbf{Sub}(G)$?\index{amenable}
\end{PROB}

The answer is positive for a large number of natural examples, see~\cite{CaMo_RelAmen} (which contains a detailed discussion of the problem) and~\cite{Wesolek_RelAmen}. A theorem of P.~Wesolek ensures moreover that if $G$ is second countable, then the set of closed amenable subgroups of $G$ is Borel, see~\cite[The\-o\-rem~A.1]{BaDuLe}. We don't know the answer to Problem~\ref{prob:LimitAmenable} in the case of the Neretin group\index{Neretin group} (see Chapter~8 in the book~\cite{book}).

A closely related problem is the following. 

\begin{PROB}
	Let $G$ be a locally compact group. Is the collection of closed locally elliptic subgroups of $G$ a closed subset of the Chabauty space $\mathbf{Sub}(G)$?\index{elliptic!locally ---}\index{locally!elliptic}
\end{PROB}

The answer is positive for the Neretin group (as a consequence of~\cite[Corol\-lary~3.6]{LBW}), but open in general. For totally disconnected groups, it is equivalent to the following. 

\begin{PROB}[C. Rosendal]
	Let $G$ be a \tdlc group and $n >0$ be an integer. Is the set $P_n(G)$ consisting of the $n$-tuples of elements of $G$ contained in a common compact subgroup, closed in the Cartesian product $G^n$? 
\end{PROB}

For $n=1$, the answer is positive by a theorem of G. Willis recalled in Chapter~9 in the book~\cite{book}.

\section{$p$-adic Lie groups}\index{Lie group!p-adic@{$p$-adic ---}}

Every $p$-adic Lie group has a continuous finite-dimen\-sion\-al linear representation over $\QQ_p$ given by its adjoint action on its Lie algebra. Hence a topologically simple $p$-adic Lie group is either linear or has a trivial adjoint representation.  The linear  topologically simple $p$-adic Lie groups are classified: they are all simple algebraic groups over $\QQ_p$, and in particular compactly generated, see~\cite[Propo\-si\-tion~6.5]{CCLTV}. We do not know whether other simple $p$-adic Lie groups exist: 

\begin{PROB}
Is there a topologically simple $p$-adic Lie group whose adjoint representation is trivial? Can it be one-dimen\-sion\-al?
\end{PROB}

A $p$-adic Lie group whose adjoint representation is trivial has an abelian Lie algebra, and is thus locally abelian. In particular it is elementary of rank~$2$ (in the sense of Wesolek), see Chapter~16 in the book~\cite{book}.

\section{Profinite groups}

Thanks to the major advances due to N.~Nikolov and D.~Segal reviewed in Chapter~5 in the book~\cite{book}, the abstract algebraic structure of finitely generated profinite groups is now well understood. To what extent is the assumption of finite generation necessary in their theory? We mention some specific questions that could guide research in this direction. 

\begin{PROB}[J. Wilson]
	Can a pro-$p$ group have a non-trivial abstract quotient that is perfect?
\end{PROB}

\begin{PROB}
	Can a hereditarily just-infinite profinite group have a proper dense normal subgroup?\index{normal subgroup!dense ---}
\end{PROB}

More problems on profinite groups are included and discussed in Chapter~3 in the book~\cite{book}. 

\section{Contraction groups}

Contraction groups\index{contraction group} appear naturally in the structure theory of \tdlc groups, in the presence of automorphisms whose scale is greater than one, see Chapter~10 in the book~\cite{book}. Moreover, when the contraction group of an automorphism is closed, it is subjected to the far-reaching results from~\cite{GlocknerWillis}. However, the following basic question remains open. 

\begin{PROB}
	Let $G$ be a \tdlc group and $\alpha \in \Aut(G)$ be a contracting automorphism, i.e. an automorphism such that $\lim_{n\to \infty} \alpha^n(g) = 1$ for all $g \in G$. Assume that the exponent of $G$ is a prime power. Does it follow that $G$ is nilpotent?
\end{PROB}

The results from~\cite{GlocknerWillis} ensure that $G$ is solvable. Moreover, the problem is known to have a positive answer if $G$ is a Lie group over a local field of arbitrary characteristic, see~\cite[Application~9.2]{Glockner}.

\section{Compactly generated simple groups}

Let $\mathscr{S}$ denote the class of non-discrete compactly generated \tdlc groups that are topologically simple. 
One of the recent trends in the structure theory of \tdlc groups is the approach that relates  the global properties of compactly generated simple \tdlc groups with the structural properties of its compact open subgroups. This was initiated in~\cite{Willis07} and~\cite{BEW} and further elaborated in~\cite{CRW} (see Chapters~17 and~18 in the book~\cite{book}). Generally speaking, a property verified by all sufficiently small compact open subgroups is called \textbf{local}. One is thus interested in relating the local and global structures of groups in $\mathscr{S}$. 

A  basic question is to evaluate the number of local isomorphism classes: two groups are called \textbf{locally isomorphic} if they contain isomorphic open subgroups. 

\begin{PROB}
	Is the number of local isomorphism classes of groups in $\mathscr{S}$ uncountable? 
\end{PROB}

An equivalent way to think of the local approach to the study of the class $\mathscr{S}$ is to ask which profinite groups embed as a compact open subgroup in a group in $\mathscr{S}$. Despite important recent progress,  our understanding of that problem remains elusive, as   illustrated by the following. 

\begin{PROB}[Y. Barnea, M. Ershov] \label{pb:BarneaErshov}
	Can a group in $\mathscr{S}$ have a compact open subgroup isomorphic to a  free pro-$p$ group? 
\end{PROB}

A basic observation from~\cite{CRW} ensures that it cannot be a free profinite group, since every group in $\mathscr{S}$ is locally pro-$\pi$ for a finite set of primes $\pi$.

Another natural problem occurring in the realm of simple \tdlc groups is the difference between topological simplicity (every closed normal subgroup is trivial) and abstract simplicity (every normal  subgroup of the underlying abstract group is trivial). Examples show that a topologically simple \tdlc group can fail to be abstractly simple\index{simple!abstractly ---}, see~\cite{Willis07} or the introduction of~\cite{CRW}. However, no compactly generated such example is known.

\begin{PROB} 
	Can a group in $\mathscr{S}$ have a proper dense normal subgroup?\index{normal subgroup!dense ---}
\end{PROB}

The following problem is closely related, see~\cite{CRW} and~\cite[Appendix B]{CaMo}. 

\begin{PROB} 
	Can a group $G$ in $\mathscr{S}$ be such that $\Inn(G)$ is not closed in $\Aut(G)$?
\end{PROB}

It has been proved in~\cite{CRW_TitsCore}   that a topologically simple \tdlc group having an element with a non-trivial contraction group has a smallest abstract normal subgroup, which is moreover simple. It is thus desirable to know whether all groups in $\mathscr{S}$ admit such an element. That property could fail if a group in $\mathscr{S}$ satisfied either of the following two conditions.  

\begin{PROB} 
	Can a group in $\mathscr{S}$ have a dense conjugacy class? Can all its closed subgroups be unimodular?
\end{PROB}

We believe that a better understanding of the class $\mathscr{S}$ as a whole necessitates to develop an intuition based on a larger pool of examples. The results from~\cite{CRW} show that all groups in $\mathscr{S}$ with a non-trivial centraliser lattice share some fundamental features with the full automorphism group of a tree. On the other hand, as long as the centraliser lattice is trivial, many of the tools developed in loc. cit. become inefficient or even useless.  

\begin{PROB}
	Find new examples of groups in $\mathscr{S}$ whose centraliser lattice is trivial. 
\end{PROB}

Another intriguing direction to explore is the relation between the algebraic structure of a group in $\mathscr{S}$ and its analytic properties, and in particular its unitary representations. In the classical case of simple Lie and algebraic groups, Kazhdan's property (T) is a landmark that is also a gateway to numerous fascinating rigidity phenomena. Some Kac--Moody groups in $\mathscr{S}$ also enjoy property (T). We do not know whether a group in $\mathscr{S}$ with a non-trivial centraliser lattice can have (T). 

\begin{PROB}
	Find new examples of groups in $\mathscr{S}$ satisfying Kazhdan's property (T).\index{Kazhdan property}
\end{PROB}

One of the main results from~\cite{CRW} is that a group in $\mathscr{S}$ with a non-trivial centraliser lattice is not amenable. The question of the existence of an infinite finitely generated simple amenable group was solved positively in~\cite{JuMo}. Its non-discrete counterpart remains open. 

\begin{PROB} \label{pb:AmenableSimple}
	Can a group in $\mathscr{S}$ be amenable?\index{amenable}
\end{PROB}

If this question has a negative answer, then none of the topological full groups considered in~\cite{JuMo} admits proper, infinite commensurated subgroups. Additionally, no non-virtually abelian finitely generated just infinite amenable group would appear as a lattice in a non-trivial way.

Amenable groups in $\mathscr{S}$ would be highly interesting since their behaviour would necessarily be very different from that of all currently known examples. On the other hand, a negative answer to the previous problem would have far-reaching consequences on discrete amenable groups. 

For a more detailed discussion of the class $\mathscr{S}$ and open problems about it, we refer to~\cite{CapECM}.

\section{Lattices}

A fundamental impetus to the study of simple \tdlc groups beyond the case of algebraic groups was the ground-breaking work of M.~Burger and S.~Mozes on lattices in products of trees, see~\cite{BuMo2} and Chapters~6 and~12 in the book~\cite{book}.  However, the mechanisms responsible for the existence or  non-existence of lattices in general simple groups remain largely mysterious\index{lattice!non-existence thereof}. The following vague problem consists in investigating that question. 

\begin{PROB}
	Which groups in $\mathscr{S}$ contain lattices? Which products of groups in $\mathscr{S}$ contain lattices with dense projections? For a group (or a product of groups) in $\mathscr{S}$, can one classify its lattices up to commensurability?
\end{PROB}

Even in the basic case (considered by Burger--Mozes) of a product  $G = G_1 \times G_2$ of two groups in $\mathscr{S}$, each acting properly and cocompactly on a regular tree, it is not clear how the existence of an irreducible cocompact lattice is reflected by the global structure of $G$ (see the Basic Question on p.~5 in \cite{BurgerMozesZimmer}). 

Another fundamental problem, also naturally suggested by the work of Burger and Mozes, is the following.

\begin{PROB}\label{pb:MultipleProd}
Let $G = G_1 \times \dots \times G_n$ be the product of $n$ non-discrete compactly generated topologically simple locally compact groups and $\Gamma \leq G$ be a lattice. Assume that the projection of $\Gamma$ to every proper subproduct of $G$ has dense image.  Can $\Gamma$ be non-arithmetic and residually finite if $n \geq 2$? Can $\Gamma$ be   simple if $n \geq 3$? 
\end{PROB}	
	
Partial results when the factors are certain Kac--Moody groups in $\mathscr{S}$ have been  established in~\cite{CaMoKM}. The special case of Problem~\ref{pb:MultipleProd} where each $G_i$ is a closed subgroup of the automorphism group of a locally finite tree $T_i$ acting cocompactly (and even $2$-transitively at infinity) is already highly interesting. 

A related problem consists in finding `exotic' lattices in the full automorphism group of a simple algebraic group over a local field. 

\begin{PROB}[{\cite[Annexe A, Problem 1]{CorHab}}]
	Let $\mathbf G$ be a simply connected absolutely simple  algebraic group over a local field $k$, of $k$-rank~$\geq 2$. Has every lattice finite image in $\mathrm{Out}(\mathbf G(k))$?
\end{PROB}

In additive combinatorics, there is a long history of considering sets that are just-not-quite groups; the modern notion of \textbf{approximate groups} was introduced by T.~Tao in~\cite{Tao_prod}. In a recent preprint~\cite{BjoHar}, M.~Bj\"orklund and T.~Hartnick consider certain subsets of a locally compact group $G$ which they call \textbf{uniform approximate lattices}. They further investigate three different tentative definitions of (non-uniform) ``approximate lattices''.

\begin{PROB}[Bj\"orklund--Hartnick]
Which is the ``right'' definition of approximate lattices?
\end{PROB}	

More specifically, we could ask for a definition such that (i)~uniform approximate lattices are approximate lattices; (ii)~a subgroup of $G$ is an approximate lattice if and only if it is a lattice.


\section{Commensurated subgroups and commensurators}

The class of \tdlc groups is closely related to the class of pairs $(\Gamma,  \Lambda)$ consisting of a discrete group $\Gamma$ and a commensurated subgroup $\Lambda \leq \Gamma$, see for example~\cite[Section~3]{ShaWil}. 
	
\begin{PROB}[Margulis--Zimmer conjecture]
Let $G$ be a connected semisimple Lie group with finite centre and $\Gamma \leq G$ be an irreducible lattice. Assume that the rank of $G$ is at least~$2$. Prove that every commensurated subgroup of $\Gamma$ is either finite or of finite index.
\end{PROB}

See~\cite{ShaWil} for an extended  discussion and partial   results. We emphasise that the only known cases where that problem has been solved concern  non-uniform lattices: there is not a single example of a uniform lattice for which the Margulis--Zimmer conjecture has been proved. 

Instead of looking for commensurated subgroups in a given group, one can dually consider the largest group in which a given group $G$ embeds as a commensurated subgroup. That group is called the \textbf{group of abstract commensurators} of $G$. It is defined in Section~6, Appendix~B of~\cite{BaK}, where the idea of the concept is attributed to J.-P. Serre and W. Neumann. 

\begin{PROB}[A. Lubotzky] \label{prob:Lub}
	Let $F$ be a non-abelian free group of finite rank. Is the group  of abstract commensurators of $F$ a simple group? 
\end{PROB}

That group of abstract commensurators is countable, but not finitely generated, see~\cite{BartBogo}. 
Several variants of that problem can be envisioned. The following one is due to  A.~Lubotzky,  S.~Mozes and R. Zimmer (Remark~2.12(i) in~\cite{LMZ}): is the relative commensurator of $F$ in the full automorphism group of its Cayley tree virtually simple? We refer to the appendix of \cite{Radu} for partial answers to the latter question as well as Problem~\ref{prob:Lub}.    It is also natural to ask whether the group of abstract commensurators of the profinite (resp. pro-$p$) completion of $F$ is simple, or whether it is topologically simple with respect to the natural \tdlc group topology that it carries (see~\cite{BEW}). The latter question is thus related to Problem~\ref{pb:BarneaErshov}.

\section{Unitary representations and  C*-simplicity}

The problems in this section pertain to a general research direction  consisting in relating the intrinsic algebraic/geometric/dynamical structure of a locally compact group with the properties of its unitary representations. Given the difficulty and depth of the theory of unitary representations of semi-simple groups over local fields, it is of course not realistic at this stage to hope for a meaningful general theory. However, recent breakthroughs suggest that some specific questions could be solved. We mention a few of them.

\begin{PROB}
	Let $G$ be a \tdlc group. Characterise the C*-simplicity of $G$ in terms of its Furstenberg boundary. 
\end{PROB}

For discrete groups, such a characterisation has been obtained recently by M.~Kalantar and M.~Kennedy: they proved in~\cite{KalKen} that a discrete group is C*-simple if and only if its action on its Furstenberg boundary is topologically free.

\medskip

The following problem was suggested by T. Steger, who reported that C.~Nebbia asked it in the 1990s. We recall that a representation of a \tdlc group is called \textbf{admissible} if the subspace of fixed points of every compact open subgroup is finite-dimen\-sion\-al.

\begin{PROB}[C. Nebbia] \label{prob:Nebbia}
	Let $T$ be a locally finite leafless tree  and $G \leq \Aut(T)$ be a closed subgroup acting $2$-transitively on the set of ends $\partial T$. Is every continuous irreducible  unitary representation of $G$ admissible?
\end{PROB}

A classical criterion (see~\cite{Dixmier} or Thm.~2.2 in~\cite{Cio}) implies that a \tdlc group all of whose  continuous irreducible unitary representations are admissible, is of type I, i.e. all of its continuous unitary representations generate a von Neumann algebra of type I. Problem~\ref{prob:Nebbia} thus leads us naturally to the following. 

\begin{PROB}  \label{prob:TypeI}
	Let $T$ be a locally finite leafless tree. Is it true that  $G \leq \Aut(T)$ is of type~I if and only if $G$ is $2$-transitive on the set of ends?
\end{PROB}

Thus a positive solution to Problem~\ref{prob:Nebbia} implies that the `if' part of Problem~\ref{prob:TypeI} holds. The converse implication in Problem~\ref{prob:TypeI} was recently proved by C.~Houdayer and S.~Raum (see~\cite{HouRau}, which also contains an extensive discussion of Problem~\ref{prob:TypeI}). 

\medskip

A topological group is called \textbf{unitarisable} if all its uniformly bounded continuous representations on Hilbert spaces are conjugate to unitary representations. Following work of B.~Sz.-Nagy, it was observed in 1950 by M.~Day,  J.~Dixmier, M.~Nakamura and Z.~Takeda that this property holds for amenable groups (their argument for discrete groups holds unchanged for topological groups). These authors raised the question whether, conversely, only amenable groups are unitarisable. This question, surveyed in~\cite{Pisier}, has led to deep results by G.~Pisier. Despite modest contributions by other authors in more recent years~\cite{EpsMon},~\cite{MonOza}, it remains completely open.

The unitarisability question makes sense more generally for locally compact groups, see~\cite{GheMon} for a partial result (beyond the locally compact setting, there are amenable topological groups without any uniformly bounded continuous representation~\cite{Ghe}). However, even the most basic tools to study it appear to fall short in the non-discrete setting. For instance:

\begin{PROB}
	Let $G$ be a unitarisable \tdlc group and $H<G$ a closed subgroup. Must $H$ also be unitarisable?
\end{PROB}

We do not even know the following particular cases:

\begin{PROB}
	Can a unitarisable \tdlc group contain a discrete non-abelian free subgroup? Can it contain a non-abelian free subgroup as a lattice?
\end{PROB}


\section{Elementary groups}

 When trying to decompose a general locally compact group into `atomic building blocks' by means of subnormal series, several families of subquotients appear to be unavoidable:   discrete groups, compact groups and topologically characteristically simple groups, see~\cite{CaMo}.  The class of elementary groups\index{elementary!group} was  introduced and studied by P.~Wesolek~\cite{Wesolek_elem} as a tool to investigate  general \tdlc groups by understanding which of them are exclusively built out of discrete and compact pieces. In that sense, those \tdlc groups are the most elementary, whence the choice of terminology. 
 
 The most general decomposition results on arbitrary 
  \tdlc groups have been obtained in the past two years by C.~Reid and P.~Wesolek in a deep and far-reaching theory which highlights the key role played by elementary groups, see~\cite{RW_chief},~\cite{RW_compression} and references therein. We thank both of them for their suggestions about the present subsection; most of the problems selected here are due to them.

 By definition, an elementary group is constructed by means of an iterative procedure involving more and more building blocks. The complexity of the resulting group is measured by an ordinal-valued  function called the \textbf{rank}, see~\cite{Wesolek_elem}.

 \begin{PROB}
 	Do there exist second countable elementary groups of arbitrarily large rank below $\omega_1$?
 \end{PROB}

 A basic fact is that a group in $\mathscr{S}$ is not elementary. Moreover, since the class of elementary groups is stable under passing to closed subgroups and Hausdorff quotients, it follows that an elementary group cannot have a subquotient isomorphic to a group in $\mathscr{S}$. We do not know whether the converse holds.

\begin{PROB}\label{pb:SimpleSubquotient}
Let $G$ be a second countable \tdlc group that is not elementary. Must there exist closed subgroup $H, K$ of $G$ with $K$ normal in $H$ such that $H/K$ belongs to $\mathscr{S}$?   
\end{PROB}

Another fundamental problem is to understand which elementary groups are topologically simple. The known simple elementary groups include all discrete simple groups, as well as various topologically simple locally elliptic groups. All of them are thus of rank~$\leq 2$. 

\begin{PROB}
	Characterise the elementary groups that are topologically simple. Are they all of rank~$2$?
\end{PROB}

The results from~\cite{CRW} show that  for many groups in $\mathscr{S}$, the conjugation action of the group on its closed subgroups has interesting dynamics. To what extent is that feature shared by non-elementary groups? The following problems are guided by this vague question.

\begin{PROB} 
	Characterise the \tdlc groups all of whose closed subgroups are unimodular. Are they all elementary? 
\end{PROB}



\begin{PROB}
Characterise the \tdlc group $G$ such that the only minimal closed invariant subsets of the Chabauty space $\mathbf{Sub}(G)$ are $\{1\}$ and $\{G\}$. Are they all elementary? 
\end{PROB}	

We refer to~\cite{LBMB} for very recent exemples.


Finally, a question of P. Wesolek blending the notion of elementarity discussed here with the classical notion of elementary amenability (of discrete groups) is as follows.


\begin{PROB}[P. Wesolek]\label{pb:AE}
 	Let $G$ be an amenable second countable \tdlc group. Must $G$ be elementary? 
\end{PROB}

Notice that a positive answer to that question would imply a negative answer to Problem~\ref{pb:AmenableSimple}. Moreover, in case the answer to Problem~\ref{pb:SimpleSubquotient} is positive, then Problems~\ref{pb:AmenableSimple} and~\ref{pb:AE} are then formally equivalent.

A more specific sub-question of Problem~\ref{pb:AE} is:

\begin{PROB}[P. Wesolek] \label{pb:Subexponential}
	Let $G$ be a compactly generated \tdlc group of subexponential growth. Must $G$ be elementary?
\end{PROB}

We refer to \cite{Cornulier:wreath} for recent examples of \tdlc groups of subexponential growth that are elementary, but not compact-by-discrete. 


\section{Galois groups}\index{Galois group}

Let $K/k$ be a field extension and $G = \Aut(K/k)$ its automorphism group, i.e. the set of those automorphisms of $K$ acting trivially on $k$. Endow $G$ with the topology of pointwise convergence. It is well known that if $K/k$ is algebraic, then $G$ is a profinite group. More generally, if $K/k$ is of finite transcendence degree, then $G$ is a \tdlc group, which is moreover discrete if and only if   $K$ is finitely generated over $k$   (see Chapter~6, \S6.3 in~\cite{Shimura}).  Thus infinitely generated transcendental field extensions of finite transcendence degree provide a natural source of examples of non-discrete \tdlc groups. While it is known that every profinite group is the Galois group of some algebraic extension (see~\cite{Waterhouse}), the corresponding problem does not seem to have been addressed for non-compact groups.

\begin{PROB}
Which \tdlc groups are Galois groups? 	
\end{PROB}

Some natural field extensions moreover yield simple groups. This is for example the case if the   fields $k$ and $K$ are algebraically closed of characteristic~$0$ 
as soon as $K/k$ is non-trivial (hence transcendental): Indeed, by~\cite[The\-o\-rem~2.9]{Rovinsky} the Galois group  $\Aut(K/k)$ has a  characteristic open subgroup which is topologically simple.

\begin{PROB}
Which topologically simple \tdlc groups are Galois groups? Which topologically simple \tdlc groups continuously embed in Galois groups?	
\end{PROB}



\begin{thebibliography}{99}
\bibitem{BaDuLe}
{\sc Bader, U., Duchesne, B., and L\'ecureux, J.}
\newblock Amenable invariant random subgroups.
\newblock {\em Israel J. Math. 213}, 1 (2016), 399--422.
\newblock With an appendix by Phillip Wesolek.

\bibitem{BaK}
{\sc Bass, H., and Kulkarni, R.}
\newblock Uniform tree lattices.
\newblock  {\em J. Amer. Math. Soc. 3}, 4 (1990),  843--902.

\bibitem{BEW}
{\sc Barnea, Y., Ershov, M., and Weigel, T.}
\newblock Abstract commensurators of profinite groups.
\newblock {\em Trans. Amer. Math. Soc. 363}, 10 (2011), 5381--5417.

\bibitem{BartBogo}
{\sc Bartholdi, L., and Bogopolski, O.}
\newblock On abstract commensurators of groups.
\newblock {\em J. Group Theory 13}, 6 (2010), 903--922.

\bibitem{BjoHar}
{\sc Bj\"orklund, M., and Hartnick, T.}
\newblock Approximate lattices.
\newblock Preprint, arXiv:1612.09246, 2016.

\bibitem{BuMo2}
{\sc Burger, M., and Mozes, S.}
\newblock Lattices in product of trees.
\newblock {\em Inst. Hautes \'Etudes Sci. Publ. Math.} 92 (2000), 151--194
  (2001).
  
\bibitem{BurgerMozesZimmer}
{\sc Burger, M.,   Mozes, S. and Zimmer, R. J.}
\newblock Linear representations and arithmeticity of lattices in products of trees. 
\newblock {\em Essays in geometric group theory}, Ramanujan Math. Soc. Lect. Notes Ser., 9 (2009),  1--25. 

\bibitem{CapECM}
{\sc Caprace, P.-E.}
\newblock Non-discrete simple locally compact groups.
\newblock Preprint, to appear in the Proceedings of the 7th European Congress
  of Mathematics, 2016.

\bibitem{CaMo}
{\sc Caprace, P.-E., and Monod, N.}
\newblock Decomposing locally compact groups into simple pieces.
\newblock {\em Math. Proc. Cambridge Philos. Soc. 150}, 1 (2011), 97--128.

\bibitem{CaMoKM}
{\sc Caprace, P.-E., and Monod, N.}
\newblock A lattice in more than two {K}ac-{M}oody groups is arithmetic.
\newblock {\em Israel J. Math. 190\/} (2012), 413--444.

\bibitem{CaMo_RelAmen}
{\sc Caprace, P.-E., and Monod, N.}
\newblock Relative amenability.
\newblock {\em Groups Geom. Dyn. 8}, 3 (2014), 747--774.

\bibitem{book}
{\sc Caprace, P.-E., and Monod, N.} (eds.)
\newblock {\em New Directions in Locally Compact Groups}.
\newblock London Mathematical Society Lecture Note Series, vol.~447
\newblock Cambridge University Press, Cambridge, 2018.

\bibitem{CRW_TitsCore}
{\sc Caprace, P.-E., Reid, C.~D., and Willis, G.~A.}
\newblock Limits of contraction groups and the {T}its core.
\newblock {\em J. Lie Theory 24}, 4 (2014), 957--967.

\bibitem{CRW}
{\sc Caprace, P.-E., Reid, C.~D., and Willis, G.~A.}
\newblock Locally normal subgroups of totally disconnected groups. {P}art {II}:
  {C}ompactly generated simple groups.
\newblock {\em Forum Math. Sigma 5\/} (2017), e12, 89.

\bibitem{Cio}
{\sc Ciobotaru, C.}
\newblock A note on type {I} groups acting on d-regular trees.
\newblock Preprint, arXiv:1506.02950, 2015.

\bibitem{CCLTV}
{\sc Cluckers, R., Cornulier, Y., Louvet, N., Tessera, R., and Valette, A.}
\newblock The {H}owe-{M}oore property for real and {$p$}-adic groups.
\newblock {\em Math. Scand. 109}, 2 (2011), 201--224.

\bibitem{CorHab}
{\sc Cornulier, Y.}
\newblock Aspects de la g\'eom\'etrie des groupes.
\newblock M\'emoire d'habilitation \`a diriger des recherches, Universit\'e
  Paris-Sud 11, 2014.

\bibitem{Cornulier:wreath}
{\sc Cornulier, Y.}
\newblock Locally compact wreath products. 
\newblock Preprint , arXiv:1703.08880, 2017.

\bibitem{Dixmier}
{\sc Dixmier, J.}
\newblock {\em {$C\sp*$}-algebras}.
\newblock North-Holland Publishing Co., Amsterdam-New York-Oxford, 1977.
\newblock Translated from the French by Francis Jellett, North-Holland
  Mathematical Library, Vol. 15.

\bibitem{EpsMon}
{\sc Epstein, I., and Monod, N.}
\newblock Nonunitarizable representations and random forests.
\newblock {\em Int. Math. Res. Not. IMRN} 22 (2009), 4336--4353.

\bibitem{Ghe}
{\sc Gheysens, M.}
\newblock Inducing representations against all odds.
\newblock Thesis (Ph.D.)--EPFL, 2017.

\bibitem{GheMon}
{\sc Gheysens, M., and Monod, N.}
\newblock Fixed points for bounded orbits in {H}ilbert spaces.
\newblock {\em Ann. Sci. \'Ec. Norm. Sup\'er. (4) 50}, 1 (2017), 131--156.

\bibitem{Glockner}
{\sc Gl\"ockner, H.}
\newblock Invariant manifolds for analytic dynamical systems over ultrametric
  fields.
\newblock {\em Expo. Math. 31}, 2 (2013), 116--150.

\bibitem{GlocknerWillis}
{\sc Gl\"ockner, H., and Willis, G.~A.}
\newblock Classification of the simple factors appearing in composition series
  of totally disconnected contraction groups.
\newblock {\em J. Reine Angew. Math. 643\/} (2010), 141--169.

\bibitem{HouRau}
{\sc Houdayer, C., and Raum, S.}
\newblock Locally compact groups acting on trees, the type {I} conjecture and
  non-amenable von {N}eumann algebras.
\newblock Preprint, arXiv:1610.00884, 2016.

\bibitem{JuMo}
{\sc Juschenko, K., and Monod, N.}
\newblock Cantor systems, piecewise translations and simple amenable groups.
\newblock {\em Ann. of Math. (2) 178}, 2 (2013), 775--787.

\bibitem{KalKen}
{\sc Kalantar, M., and Kennedy, M.}
\newblock Boundaries of reduced {$C^*$}-algebras of discrete groups.
\newblock {\em J. Reine Angew. Math. 727\/} (2017), 247--267.

\bibitem{LBMB}
{\sc Le~Boudec, A., and Matte Bon, N.}
\newblock Locally compact groups whose ergodic or minimal actions are all free
\newblock Preprint, arXiv:1709.06733, 2017.

\bibitem{LBW}
{\sc Le~Boudec, A., and Wesolek, P.}
\newblock Commensurated subgroups in tree almost automorphism groups.
\newblock Preprint, arXiv:1604.04162, 2016.

\bibitem{LMZ}
{\sc Lubotzky, A., Mozes, S., and Zimmer, R.}
\newblock  Superrigidity for the commensurability group of tree lattices.
\newblock {\em Comment. Math. Helv. 69}, 4 (1994),  523--548. 


\bibitem{MonOza}
{\sc Monod, N., and Ozawa, N.}
\newblock The {D}ixmier problem, lamplighters and {B}urnside groups.
\newblock {\em J. Funct. Anal. 258}, 1 (2010), 255--259.

\bibitem{Pisier}
{\sc Pisier, G.}
\newblock {\em Similarity problems and completely bounded maps}, expanded~ed.,
  vol.~1618 of {\em Lecture Notes in Mathematics}.
\newblock Springer-Verlag, Berlin, 2001.
\newblock Includes the solution to ``The Halmos problem''.

\bibitem{Radu}
{\sc Radu, N.}
\newblock {New simple lattices in products of trees and their projections}.
\newblock ArXiv preprint 1712.01091, 2017. 

 
\bibitem{RW_chief}
{\sc {Reid}, C.~D., and Wesolek, P.~R.}
\newblock The essentially chief series of a compactly generated locally compact
  group.
\newblock ArXiv preprint 1509.06593, 2015.

\bibitem{RW_compression}
{\sc {Reid}, C.~D., and Wesolek, P.~R.}
\newblock Dense normal subgroups and chief factors in locally compact groups.
\newblock ArXiv preprint 1601.07317, 2016.

\bibitem{Rovinsky}
{\sc Rovinsky, M.}
\newblock Motives and admissible representations of automorphism groups of
  fields.
\newblock {\em Math. Z. 249}, 1 (2005), 163--221.

\bibitem{ShaWil}
{\sc Shalom, Y., and Willis, G.~A.}
\newblock Commensurated subgroups of arithmetic groups, totally disconnected
  groups and adelic rigidity.
\newblock {\em Geom. Funct. Anal. 23}, 5 (2013), 1631--1683.

\bibitem{Shimura}
{\sc Shimura, G.}
\newblock {\em Introduction to the arithmetic theory of automorphic functions}.
\newblock Publications of the Mathematical Society of Japan, No. 11. Iwanami
  Shoten, Publishers, Tokyo; Princeton University Press, Princeton, N.J., 1971.
\newblock Kan\^o Memorial Lectures, No. 1.

\bibitem{Tao_prod}
{\sc Tao, T.}
\newblock Product set estimates for non-commutative groups.
\newblock {\em Combinatorica 28}, 5 (2008), 547--594.

\bibitem{Waterhouse}
{\sc Waterhouse, W.~C.}
\newblock Profinite groups are {G}alois groups.
\newblock {\em Proc. Amer. Math. Soc. 42\/} (1973), 639--640.

\bibitem{Wesolek_elem}
{\sc Wesolek, P.}
\newblock Elementary totally disconnected locally compact groups.
\newblock {\em Proc. Lond. Math. Soc. (3) 110}, 6 (2015), 1387--1434.

\bibitem{Wesolek_RelAmen}
{\sc Wesolek, P.}
\newblock A note on relative amenability.
\newblock {\em Groups Geom. Dyn. 11}, 1 (2017), 95--104.

\bibitem{Willis07}
{\sc Willis, G.~A.}
\newblock Compact open subgroups in simple totally disconnected groups.
\newblock {\em J. Algebra 312}, 1 (2007), 405--417.
	
	
\end{thebibliography}
\end{document}